\documentclass[12pt]{article}

\usepackage{url}

\makeatletter
\newfont{\footsc}{cmcsc10 at 8truept}
\newfont{\footbf}{cmbx10 at 8truept}
\newfont{\footrm}{cmr10 at 10truept}
\makeatother
\usepackage{amsmath}
\usepackage{amssymb}
\usepackage{exscale}
\usepackage{amsthm}
\usepackage{epsfig}
\usepackage{verbatim}
\usepackage{setspace}

\usepackage{graphics,color}

\usepackage[paper=a4paper,left=30mm,right=20mm,top=25mm,bottom=30mm]{geometry}
    \newenvironment{dedication}
        {\vspace{6ex}\begin{quotation}\begin{center}\begin{em}}
        {\par\end{em}\end{center}\end{quotation}}

\theoremstyle{plain}
\newtheorem{theorem}{Theorem}[section]
\newtheorem{proposition}[theorem]{Proposition}
\newtheorem{lemma}[theorem]{Lemma}

\theoremstyle{definition}
\newtheorem{definition}[theorem]{Definition}

\newtheorem{observation}[theorem]{Observation}
\newtheorem{remark}[theorem]{Remark}

\DeclareMathOperator{\1}{\bf{1}}

\newcommand{\N}{{\mathbb{N}}}

\newcommand{\R}{{\mathbb{R}}}

\newcommand{\PP}{{\mathbb{P}}}
\newcommand{\CC}{{\mathcal{C}}}
\newcommand{\DD}{{\mathcal{D}}}
\newcommand{\EE}{{\mathcal{E}}}

\newcommand{\X}{{\bf X}}
\newcommand{\Y}{{\bf Y}}

\def\bssigma{\boldsymbol{\sigma}}
\def\bsx{\boldsymbol{x}}
\def\bsy{\boldsymbol{y}}

\def\bsp{\boldsymbol{p}}
\def\bs0{\bf 0}
\def\eps{\varepsilon}

\newcommand{\change}[1]{{#1}}

\title{On Negative Dependence Properties of Latin Hypercube Samples and Scrambled Nets}

\author{Benjamin Doerr\thanks{Laboratoire d'Informatique (LIX), CNRS, \'Ecole Polytechnique, Institut Polytechnique de Paris, Palaiseau, France ({\tt lastname@lix.polytechnique.fr}).}
\and Michael Gnewuch\thanks{Institut f\"ur Mathematik, Universit\"at Osnabr\"uck, 
Germany ({\tt michael.gnewuch@uni-osnabrueck.de}).}}

\begin{document}

\maketitle


\begin{dedication}
Dedicated to the $60$-th birthday of Klaus Ritter
\end{dedication}

\begin{abstract}
We study the notion of $\gamma$-negative dependence of random variables. This notion is a relaxation of the notion of negative orthant dependence (which corresponds to $1$-negative dependence), but nevertheless it still ensures concentration of measure and allows to use large deviation bounds of Chernoff-Hoeffding- or Bernstein-type. 
We study random variables based on random points $P$. These random variables appear naturally in the analysis of the discrepancy of $P$ or, equivalently, of a suitable worst-case integration error of the quasi-Monte Carlo cubature that uses the points in $P$ as integration nodes.  
We introduce the correlation number, which is the smallest possible value of $\gamma$ that ensures $\gamma$-negative dependence.
We prove that the random variables of interest based on Latin hypercube sampling or on $(t,m,d)$-nets do, in general, not have a correlation number of $1$, i.e., they are not negative orthant dependent. But it is known that the random variables based on Latin hypercube sampling in dimension $d$ are actually $\gamma$-negatively dependent with $\gamma \le e^d$, and the resulting probabilistic discrepancy bounds do only mildly depend on the $\gamma$-value.
 \end{abstract}

\section{Introduction}

For the solution of discrete or continuous problems, algorithms based on independent random variables are a popular choice. This is mostly due to the following two  advantages.
\begin{itemize}
\item Stochastic independence usually leads to an easy analysis of the resulting error. In particular, typically the concentration of measure phenomenon 
allows to use helpful large deviation bounds.
\item Additionally, the actual error can be estimated statistically by repeating the computation several times. 
\end{itemize}
But, depending on the specific setting, using independent random variables may also have serious disadvantages. Let us mention some examples.
\begin{itemize}
\item Most likely, the corresponding algorithms do not achieve the best possible asymptotic convergence rate.
\item The corresponding algorithms may not provide admissable solutions. Specific examples from discrete mathematics and computer science are integer linear programs with certain hard constraints -- to satisfy these hard constraints the random variables involved necessarily have to be dependent. 
\item If the independent random variables are random sample points, they exhibit a lack of (combinatorial, algebraic or geometric) structure. For instance, independent random samples typically contain (undesired) point clusters or gaps. Furthermore, they do not possess a group structure like, e.g., lattice rules or digital $(t,m,s)$-nets, which allows to use Fourier analytical methods. 
\end{itemize}
Therefore it is an important and interesting problem to find algorithms based on dependent random variables that, on the one hand side, share the advantages of algorithms based on independent random variables (which are most notably the easy error analysis and the benefits of the concentration of measure phenomenon) and, on the other hand side, avoid for the specific problem at hand their unwanted disadvantages. This problem is further motivated in much more detail in the introductions of \cite{DDG18, GH19, Lem18} for the specific setting where one wants to approximate expectations of random variables $f$ or integrals $\int_{[0,1]^d} f(\bsx)\, d\bsx$ of integrands $f:[0,1]^d \to \R$ by estimators
\begin{equation*}
\frac{1}{N} \sum^N_{n=1} f(X_n),
\end{equation*}
where $X_1, \ldots, X_N$ are random points.  

To tackle this problem and to analyze promising randomized algorithms, we need an appropriate probabilistic framework. 
A suitable framework is the one of negatively dependent (or negatively correlated) random variables. 
We focus \change{here} on the following notion 
for Bernoulli or binary random variables, i.e., 
random variables that only take
values in $\{0,1\}$.  

\begin{definition}\label{Def_Neg_gamma_Dep}
Let $\gamma \ge 1$. 
We call binary random variables 
$T_1,T_2,\ldots,T_N$ 
\emph{upper $\gamma$-negatively dependent} if  
\begin{equation}
\label{cond2'}
\PP\left( \bigcap_{j\in J} \{T_j = 1\} \right) 
\le \gamma \prod_{j\in J} \PP( T_j = 1)
\hspace{2ex}\text{for all $\emptyset \neq J \subseteq [1..N]$,}
\end{equation}
and \emph{lower $\gamma$-negatively dependent} if  
\begin{equation}
\label{cond2''}
\PP\left( \bigcap_{j\in J} \{T_j = 0\} \right) 
\le \gamma \prod_{j\in J} \PP(T_j = 0)
\hspace{2ex}\text{for all $\emptyset \neq J \subseteq [1..N]$.}
\end{equation}
We call $T_1,T_2,\ldots,T_N$ \emph{$\gamma$-negatively dependent} if
both conditions (\ref{cond2'}) and (\ref{cond2''}) are satisfied.
\end{definition}

Essentially the same notion, called \emph{$\lambda$-correlation}, was introduced earlier, but in a formally different way in \cite{PS97}. 
$1$-negative dependence is usually called \emph{negative orthant dependence}, cf. \cite{BSS82}. To allow for a parameter $\gamma \ge 1$ (or $\lambda \ge 1$, respectively) instead of insisting on $\gamma$ to be one
relaxes the condition of negative orthant dependence, but still ensures that large deviation bounds of Chernoff-Hoeffding type hold for sums of $\gamma$-negatively dependent random variables, see \cite{PS97, GH19, Doerr20bookchapter}. 

Now one may ask the question if allowing the parameter $\gamma$ to be different from $1$ is a fruitful generalization. To the best of our knowledge, parameters $\gamma$ (or $\lambda$  for the notion of $\lambda$-correlation) different from $1$ have been employed \change{successfully} for the first time in \cite{GH19} to establish favorable pre-asymptotic upper bounds for the star discrepancy of Latin hypercube samples \cite{PS_PC20}.

Let us explain this application in more detail: 
Here we deal with binary random variables $T_i$, $i=1, \ldots,N$, of the form
$T_i = \1_A(X_i)$, where $A$ is a Lebesgue-measurable subset of $[0,1)^d$ 
whose characteristic function
is denoted by $\1_A$,  and $X_1, \ldots, X_N$ are random points in $[0,1)^d$.

\begin{definition}\label{Def_neg_dep}
For $d\in \N$ we put 
$$\mathcal{C}^d_0 : = \{ [0,a) \,|\, a\in [0,1]^d\}
\hspace{3ex}\text{and}\hspace{3ex} 
\mathcal{D}^d_0 : = \{ B\setminus A \,|\, A,B \in \mathcal{C}^d_0 \}.$$
Note that for $a = 0$, we have $\emptyset = [0,a) \in \mathcal{C}^d_0$. Hence all $B \in \mathcal{C}^d_0$ are also in $\mathcal{D}^d_0$ since $B$ can be written as $B = B \setminus \emptyset$. We thus have $\mathcal{C}^d_0 \subseteq \mathcal{D}^d_0$. 

Let $\mathcal{S}\in \{\mathcal{C}^d_0, \mathcal{D}^d_0\}$. We say that 
a random point set $\X=(X_n)^N_{n=1}$ in $[0,1)^d$ is \emph{$\mathcal{S}$-$\gamma$-negatively dependent} if for each $S\in \mathcal{S}$ the random variables 
\begin{equation*}
\1_S(X_1), \ldots, \1_S(X_N) \hspace{2ex}\text{are $\gamma$-negatively dependent.}
\end{equation*} 
\end{definition}
It was shown in \cite[Theorem~3.6]{GH19} that Latin hypercube samples $\X$ are $\mathcal{D}^d_0$-$e^{d}$-negatively dependent. Similar results for other random point sets can be found in \cite{WGH20}.
Furthermore, \cite[Theorem~4.4]{GH19} established that for each $c>0$ a $\mathcal{D}^d_0$-$\gamma$-negatively dependent random point set $\X$ satisfies the star discrepancy bound
\begin{equation}\label{est:star_disc}
D^\ast_N(\X) \le c \sqrt{\frac{d}{N}}
\end{equation}
with success probability at least $1- \gamma \exp(-(1.6741 \cdot c^2 - 10.7042)\cdot d)$. (\change{The constants in the bound of the success probability have} been improved slightly in \cite{GPW20}.) Note that the success probability depends only mildly on the parameter $\gamma$, \change{hence larger values of $\gamma$ can easily be leveraged by slightly increasing the parameter $c$.}
This is, \change{e.g., visible in} \cite[Corollary~4.5]{GH19}: Let $d, N\in \N$. Estimate \eqref{est:star_disc} implies
that there exists a realization $P_{\rm MC}$ of Monte Carlo points (i.e., $\mathcal{D}^d_0$-$1$-negatively dependent random points) having star discrepancy
$$
D^\ast_N(P_{\rm MC}) \le 2.5287 \cdot \sqrt{\frac{d}{N}},
$$ 
improving in terms of the leading constant the previous works of Heinrich, Novak, Wasil\-kowski, and Wo\'{z}niakowski~\cite{HNWW01} and Aistleitner~\cite{Ais11}, and a realization $P_{\rm LHS}$ of Latin hypercube sampling (i.e., of $\mathcal{D}^d_0$-$e^{d}$-negatively dependent random points) satisfying 
$$
D^\ast_N(P_{\rm LHS}) \le 2.6442 \cdot \sqrt{\frac{d}{N}}\,.
$$ 
Note that all these upper bounds are asymptotically tight as shown by the lower bounds~\cite{Doe14,DDG18}. For the leading constants in the upper bounds, we note that the $\gamma$-values of $1$ (Monte Carlo points) and $e^d$ (Latin hypercube sampling) differ significantly (especially in high dimension $d$), but the \change{resulting} leading constants $2.5287$ and $2.6442$ differ only marginally. (\change{Both constants} have recently been improved in \cite{GPW20} 
\change{to $2.4968$ and $2.6135$, respectively}.)

Nevertheless, the question remains if Latin hypercube samples are even $\mathcal{D}^d_0$-negatively orthant dependent and therefore introducing a parameter $\gamma >1$ is unnecessary if we just use a different proof idea than the one used in the proof of \cite[Theorem~3.6]{GH19}?

Since Latin hypercube samples are actually $\mathcal{C}^d_0$-$1$-negatively dependent, see \cite[Theorem~3.6]{GH19}, several colleagues expressed their belief that Latin hypercube samples are also $\mathcal{D}^d_0$-$1$-negatively dependent.
In Proposition~\ref{prop33} we show that this is actually not the case.
This gives strong reason to believe that there are other important families of random variables that are \emph{not} negatively orthant dependent, but \emph{still} $\gamma$-negatively dependent for some $\gamma >1$ and therefore, nevertheless, fall under the concentration of measure phenomenon and satisfy large deviation bounds.


For the quantitative analysis of negatively dependent random variables we introduce the notion of the $\mathcal{S}$-correlation number.

\begin{definition}\label{Def_corr_num}
Let $\mathcal{S}\in \{\mathcal{C}^d_0, \mathcal{D}^d_0\}$, and $\X=(X_n)^N_{n=1}$ 
be a random point set in $[0,1)^d$. 
The \emph{$\mathcal{S}$-correlation number} $\gamma_{\mathcal{S}}(\X) \in [1,\infty]$ of $\X$ is defined as the smallest $\gamma$ such that $\X$ is $\mathcal{S}$-$\gamma$-negatively dependent, i.e., 
\begin{equation*}
\gamma_{\mathcal{S}}(\X) := \sup \left\{\frac{\PP (\bigcap_{j\in J} \{\1_S( X_j) =0\})}{\prod_{j\in J} \PP(\1_S(X_j) = 0)}, \frac{\PP (\bigcap_{j\in J} \{\1_S( X_j) =1\})}{\prod_{j\in J} \PP(\1_S(X_j) = 1)} \,\middle|\, S\in \mathcal{S}, \emptyset \neq J \subseteq [1..N]\right\}.
\end{equation*}
\end{definition}


In Section~\ref{LHS} we recall the definition of Latin hypercube samples and show that their $\mathcal{D}^d_0$-correlation number grows with the dimension $d$ at least as $\Omega(\sqrt{d})$, see Proposition~\ref{prop33}. In Section~\ref{NETS} we study scrambled $(t,m,d)$-nets. The random points of scrambled nets exhibit even stronger dependencies than Latin hypercube samples. Together with lattice rules, $(t,m,d)$-nets are arguably the most popular choice of sample points used in quasi-Monte Carlo methods. Since there are many different ways known how to scramble $(t,m,d)$-nets, we first introduce the notion of an abstract scrambled $(t,m,d)$-net, see Definition~\ref{defasnet}.
This notion captures the crucial properties that random scrambling should exhibit.
Indeed, Proposition~\ref{ElemProp} shows that the most important variants of scramblings directly satisfy this notion or can easily be \change{modified} to do so.
In Remark~\ref{Rem:Symmetrization} we prove that in order to obtain small correlation numbers it is indispensable that the considered scrambled $(t,m,d)$-nets are symmetrized, i.e., that their random points are exchangeable.
In Proposition~\ref{Ex:Correlation_Constant_Nets}
we demonstrate that also abstract scrambled $(t,m,d)$-nets do, in general, not have a $\mathcal{D}^d_0$-correlation number of one.

\section{Latin Hypercube Sampling}
\label{LHS}


In this section, we show that Latin hypercube samples are not $\DD_0^d$-$1$-negatively dependent. In fact, there is no constant $\gamma$ (independent of $d$) such that we have $\DD_0^d$-$\gamma$-negative dependence. This shows, in particular, that the use of a $\DD_0^d$-correlation number larger than one in the discrepancy analysis of Latin hypercube samples in~\cite{GH19} cannot be avoided (though the gap between our lower bound of $\Omega(\sqrt d)$ and the upper bound of $\exp(d)$ used in~\cite{GH19} remains large). 

In simple words, a Latin hypercube sample of $N$ points in $[0,1)^d$ is constructed by partitioning $[0,1)^d$ into $N^d$ axis-aligned cubes isomorphic to $[0,1/N)^d$, selecting randomly $N$ such cubes in a way that no two overlap in any one-dimensional projection, and then placing one point into each selected cube independently and uniformly at random. While each single point is uniformly distributed in $[0,1)^d$, the dependent choice of the small cubes leads to dependencies among the points. The precise definition of a Latin Hypercube sampling is as follows.
    
\begin{definition}
A \emph{Latin hypercube sample} (LHS) $(X_n)_{n=1}^N$ in $[0,1)^d$ is of the form
\begin{equation*}
X_{n,j} = \frac{\pi_j(n-1) + U_{n,j}}{N},
\end{equation*}
where $X_{n,j}$ denotes the $j$th coordinate of $X_n$, $\pi_j$ is a permutation
of $0,1,\ldots, N-1$, uniformly chosen at random, and $U_{n,j}$ is uniformly distributed in $[0,1)$.
The $d$ permutations $\pi_j$ and the $dN$ random variables $U_{n,j}$ are mutually independent. 
\end{definition}

The definition of Latin hypercube sampling presented above was introduced by McKay, Beckman, 
and Conover \cite{MBC79} for the design of computer experiments.
As shown by Owen \cite{Owe97b}, an estimator $\hat{I}$ of an integral $I$ based on Latin hypercube samples with $N$ points never leads to a variance greater than that of the corresponding  estimator based on $N-1$ Monte Carlo points.
As an extreme example, we note that the one-dimensional projections of Latin hypercube samples are much more evenly distributed than the one-dimensional projections of Monte Carlo points. Indeed, for a one-dimensional Latin hypercube sample of size $N$ the star discrepancy is at most $1/N$, while for a one-dimensional Monte Carlo sample of the same size the star discrepancy is of order $1/\sqrt{N}$, cf., e.g., \cite[Remark~3.2]{GH19}.

The following theorem is a result from \cite{GH19}.

\begin{theorem}\label{Thm-d-LHS}
Let $d\in \N$. 
Let $\X := (X_n)_{n=1}^N$ be a LHS in $[0,1)^{d}$.
For $a, b \in [0,1)^d$ let $A:= [0,a)$, $B:= [0,b)$, and $D:= B \setminus A$. Let 
\begin{equation*}
\gamma_{A,B} := \prod_{i=1}^{d} \delta_i \hspace{2ex}\text{with} \hspace{2ex} \delta_i := \begin{cases} 1 & \mbox{if} \quad a_i, b_i \in \{1/N, 2/N, \ldots, 1\} \quad \mbox{or} \quad a_i=0,\\ 
 e & \mbox{else}.\end{cases}
\end{equation*}
Then the  random variables 
\begin{equation}\label{LHS_neg_dep}
\1_D(X_1), \ldots, \1_D(X_N) \hspace{2ex}\text{are $\gamma_{A,B}$-negatively dependent.}
\end{equation}
In particular, 
the random points $(X_n)_{n=1}^N$ are $\mathcal{C}^d_0$-$1$-negatively dependent as well as $\mathcal{D}^d_0$-$e^d$-negatively dependent.
\end{theorem}

We now show that if Latin hypercube samples $\X$ consisting of $N=d$ points are $\mathcal{D}^d_0$-$\gamma$-negatively dependent for some $\gamma = \gamma_d$, then $\gamma$ necessarily satisfies $\gamma_2  \ge 2$  and $\liminf_{d\to \infty} \gamma_d = \infty$ for $d$ tending to $\infty$. In particular, we cannot simply replace the $\gamma_{A,B}$ in Theorem~\ref{Thm-d-LHS} by $1$ or any other constant independent of $d$.

\begin{proposition}\label{prop33}
Let $d \ge 2$ be a natural number, and let $\X := (X_n)_{n=1}^d$ be a LHS in $[0,1)^d$. Then
\begin{equation}\label{eq:claim_LHS}
\gamma_{\mathcal{D}^2_0} (\X) \ge 2
\hspace{3ex}\text{and}\hspace{3ex}
\liminf_{d\to \infty} \frac{\gamma_{\mathcal{D}^d_0} (\X)}{ \sqrt{d}} 
\ge \sqrt{\frac{2\pi}{e}}.
\end{equation}
\end{proposition}

\begin{proof}
Let $\varepsilon \in (0,1/d)$. We put 
\begin{equation*}
A = A(\varepsilon,d) := [0, \tfrac{d-1}{d})^d
\hspace{3ex}\text{and}\hspace{3ex}
B = B(\varepsilon, d) := [0, \tfrac{d-1}{d}+ \varepsilon)^d 
\end{equation*}
as well as 
$D = D(\varepsilon,d) := B\setminus A$.
We want to verify \eqref{eq:claim_LHS} by showing that
\begin{equation}\label{eq:4}
\frac{\PP(X_1,\ldots, X_d\in D)}{\prod_{i=1}^d \PP(X_i\in D)} \ge (1 - g(d,\varepsilon)) \sqrt{\frac{2\pi}{e}} \cdot \sqrt{d},
\end{equation}
where $g(d,\varepsilon)$ is such that $\lim_{\eps \to 0} g(d,\eps) = 0$, and that for $d=2$ the left hand side of \eqref{eq:4} tends to $2$ if $\varepsilon$ tends to $0$.

Let $\pi_1, \dots, \pi_d : [0..d-1] \to [0..d-1]$ be the permutations used in the construction of the Latin hypercube sample $\X = (X_i)^d_{i=1}$. Consider the event 
\[
\EE:= \{\forall i \in [0..d-1]\, \exists j(i) \in [1..d]: \pi_{j(i)}(i) = d-1\}.
\]
To have $\{X_1,\ldots, X_d\} \subset D$, it is necessary that the event $\EE$ occurs, hence 
\begin{align*}
  \PP(X_1,\ldots, X_d\in D) &= \PP(\{X_1,\ldots, X_d\in D\} \cap \EE) \\
  &= \PP(\EE) \cdot \PP (X_1, \dots, X_d \in D \,|\, \EE).
\end{align*}
Since we use only $d$ permutations to construct $\X$, the event $\EE$ is equivalent to saying that the $\pi_j^{-1}(d-1)$, $j = 1, \dots, d$, are pairwise different. 
This observations allows to easily compute the probability of $\EE$ by
\begin{equation*}
\PP(\EE) = \frac{[d\cdot(d-1)!] \cdot [(d-1)\cdot (d-1)!] \cdots [1\cdot (d-1)!] }{[d!]^d}
= \frac{d!}{d^d}.
\end{equation*}
The event $\EE$, again by elementary counting, also implies that for each $i$ there is exactly one $j$ such that $\pi_j(i) = d-1$, that is, each point of $\X$ has exactly one coordinate in $[\frac{d-1}{d},1)$. Therefore, 
\begin{equation*}
\PP (X_i \in D \,|\, \EE) = \frac{\varepsilon}{1/d} = d\varepsilon
\hspace{3ex}\text{for all $i\in [1..d]$.}
\end{equation*}
Conditioned under $\EE$ all events $\{ X_i\in D\}$ are independent, implying
\[
\PP(X_1, \ldots, X_d \in D \,|\, \EE)  = \prod^d_{i=1} \PP(X_i \in D \,|\, \EE) = (d\eps)^d.
\]
This finally gives
\begin{equation*}
\PP(X_1,\ldots, X_d\in D)
= \PP(\EE) \cdot \PP (X_1, \dots, X_d \in D \,|\, \EE) = \frac{d!}{d^d} (d\eps)^d = d! \cdot \varepsilon^d.
\end{equation*}
On the other hand we have  
$$
\lambda^d(D) = \lambda^d(B) - \lambda^d(A)
= \left( \frac{d-1}{d} + \varepsilon \right)^d -   \left( \frac{d-1}{d} \right)^d
= \varepsilon \left( \sum_{k=0}^{d-1} {d \choose k} \left( 1-\frac{1}{d} \right)^k \varepsilon^{d-1-k} \right)
$$
and thus, for all $i\in [1..d]$,
\begin{equation*}
\PP ( X_i \in D ) = \lambda^d(D) = \varepsilon \cdot d \left( 1-\frac{1}{d} \right)^{d-1} + R(\varepsilon,d),
\end{equation*}
where
$$
R(\varepsilon,d)
= \varepsilon^{2} \left( \sum_{k=0}^{d-2} {d \choose k} \left( 1-\frac{1}{d} \right)^k \varepsilon^{d-2-k} \right).
$$
This gives us 
\begin{equation}\label{quotient_1}
\frac{\PP(X_1,\ldots, X_d\in D)}{\prod_{i=1}^d \PP(X_i\in D)} = \frac{d!}{d^d} \left[ \left(1-\frac{1}{d} \right)^{d\,} \right]^{1-d} \cdot T(\varepsilon,d),
\end{equation}
where
$$T(\varepsilon, d) := \left( 1- \frac{R(\varepsilon,d)}{\lambda^d(D)} \right)^d.$$
We now want to determine the strong  asymptotic behavior of the quotient on the left hand side in equation~\eqref{quotient_1} if $d$ tends to infinity and $\varepsilon$ tends to zero (sufficiently fast). To this purpose let us first remark that by choosing $\varepsilon = \varepsilon(d)$ small enough we get $\liminf_{d\to \infty} T(\varepsilon,d) =1$.
Furthermore, it is well-known~\cite{Robbins55} that 
\begin{equation*}
\sqrt{2\pi d}\, e^{-d} < \frac{d!}{d^d} < \sqrt{2\pi d} \,e^{-d} \exp \left( \tfrac{1}{12d} \right).
\end{equation*}
Finally, we get $(1-1/d)^d \to e^{-1}$ as $d\to\infty$ and
\begin{equation*}
\begin{split}
e^{-d} \left[ \left(1-\frac{1}{d} \right)^d \right]^{-d} 
&= \exp \left( -d \left( 1+ d\ln \left( 1 - \frac{1}{d} \right) \right) \right)\\
&= \exp \left( -d \left(1+ d \left( -\frac{1}{d} - \frac{1}{2d^2} - \frac{1}{3d^3} - \ldots \right)\right)\right)\\
&= \exp \left( \frac{1}{2} + \frac{1}{3d} + \frac{1}{4d^2} + \ldots \right)\\
&= e^{1/2} \cdot \exp \left( \frac{1}{d} \left( \frac{1}{3} + \frac{1}{4d} +  \frac{1}{5d^2} + \ldots \right)\right),
\end{split}
\end{equation*}
and the last term converges to $e^{1/2}$ as $d$ tends to infinity.
Thus we get altogether 
\begin{equation}
\frac{\PP(X_1,\ldots, X_d\in D(\varepsilon, d))}{\prod_{i=1}^d \PP(X_i\in D(\varepsilon, d))} \approx \sqrt{\frac{2\pi}{e}} \cdot \sqrt{d}
\end{equation}
for $d\to \infty$ in the strong asymptotic sense as long as $\varepsilon = \varepsilon(d) \to 0$ sufficiently fast.

In the case $d=2$ our calculations above lead to $\PP(X_1, X_2 \in D) = 2 \varepsilon^2$ and 
$$\PP(X_i \in D) = \lambda^2(D) =  \left(\frac{1}{2} + \varepsilon \right)^2 - \left(\frac{1}{2} \right)^2
=\varepsilon(1 + \varepsilon)
\hspace{3ex}\text{ for $i=1,2$.}$$
Hence
$$\lim_{\varepsilon \to 0} \frac{\PP(X_1, X_2 \in D(\varepsilon, 2))}{\PP(X_1 \in D(\varepsilon, 2)) \PP(X_2 \in D(\varepsilon, 2))} 
= \lim_{\varepsilon \to 0} \frac{2}{(1+ \varepsilon)^2} = 2.$$
\end{proof}

\section{Scrambled Nets}
\label{NETS}

In this section, we show that also the dependencies stemming from scrambled net constructions can lead to a $\mathcal{D}^d_0$-correlation number larger than one. We give an example that is not just a Latin hypercube sample, and further, where the dependencies of Latin hypercube samples lead to a truly negative correlation. This indicates that the stronger dependencies of scrambled nets can be more harmful than those of Latin hypercubes.

\begin{definition}
Let $d,N\in\N$.
We call an $N$-point multi-set (or family) $P= (p^{(n)})^{N}_{n=1}$ in $[0,1)^d$ \emph{fair with 
respect to a Borel measurable subset $F$} of $[0,1)^d$ 
if $N \lambda^d(F) = |\{ n \,|\, p^{(n)} \in F\}|$. In this situation we also say that \emph{$F$ is fair 
with respect to $P$}.
\end{definition}


\begin{definition}\label{defnet}
Let us fix integers $d,b \ge 2$, and $m \ge t \ge 0$.
\begin{itemize}
\item[(i)] An \emph{elementary interval} or \emph{canonical box} $E$ in base $b$ of order $m-t$ is a subset of $[0,1)^d$ of the form
\begin{equation*}
E= [k_1 b^{-\ell_1}, (k_1+1) b^{-\ell_1}) \times \cdots \times
[k_d b^{-\ell_d}, (k_d+1) b^{-\ell_d})
\end{equation*}
for non-negative integers $\ell_i$ and $k_i < b^{\ell_i}$, $i=1,\ldots,d$,
that has volume $\lambda^d(E) = b^{-(m-t)}$.
\item[(ii)] Let $N = b^m$. An $N$-point multi-set $P= (p^{(n)})^{N}_{n=1}$ in $[0,1)^d$ 
is called a \emph{$(t,m,d)$-net in base $b$} if every elementary
interval in base $b$ of order $m-t$ is fair with respect to $P$, i.e., if it contains exactly $b^t$ 
points of $P$.
\end{itemize}
\end{definition}

The construction of nets goes back to Sobol' \cite{Sob67} and 
Faure \cite{Fau82}. The concept of $(t,m,d)$-nets was introduced by Niederreiter in
\cite{Nie87} and studied in depths by him and many other researchers, see, e.g., 
\cite{DP10, Nie92} and the literature
mentioned therein. As these works show, the property of being a $(t,m,d)$-net is sufficient to lead to a small discrepancy. In practical applications, one often also wants to have random low-discrepancy point sets, where random means that the points to a certain degree resemble uniformly distributed points. To capture the important properties of existing randomized variants of $(t,m,d)$-nets, we formulate the following definition. 

\begin{definition}\label{defasnet}
  Let $b, d, m, t$, and $N$ as in Definition~\ref{defnet}. For all $k_1, \dots, k_d \in [0..b^{m-t}-1]$, the set $C = \prod_{j=1}^d [k_j b^{-(m-t)}, (k_j+1) b^{-(m-t)})$ is called a \emph{basic cube}. Let $\Y = (Y_i)_{i=1}^N$ be a tuple of $N$ random points in $[0,1)^d$. We say that $\Y$ is an \emph{abstract scrambled $(t,m,d)$-net in base~$b$} if it satisfies the following three conditions.
  \begin{itemize}
  \item[(i)] With probability one, $(Y_n)^N_{n=1}$ is a $(t,m,d)$-net in base $b$.
  \item[(ii)] For each $i \in [1..N]$ and each basic cube $C$, we have $\PP(Y_i \in C) = b^{-d(m-t)}$.
  \item[(iii)] For all measurable $M \subseteq [0,1)^d$, all $\emptyset \neq J \subseteq [1..N]$, and all families of basic cubes $(C_j)_{j\in J}$ that satisfy 
  $\PP (\forall j \in J : Y_j \in C_j)> 0$, we have
\[\PP(\forall j \in J : Y_j \in M \mid \forall j \in J : Y_j \in C_j) = \prod_{j\in J}   \frac{\lambda^d(M \cap C_j)}{\lambda^d(C_j)}.\] 
  \end{itemize}
\end{definition}

We note that condition~(ii) and condition~(iii) immediately imply the condition 
\begin{itemize}
 \item[(iv)] For each $i \in [1..N]$, the point $Y_i$ is uniformly distributed in $[0,1)^d$.
 \end{itemize}

Note that our notion of abstract scrambled $(t,m,d)$-nets is defined by (desirable) properties of an $N$-tuple of random points $\Y$ rather than by a particular way of scrambling. Nevertheless, in practice one often starts with a (deterministic) $(t,m,d)$-net and uses than a specific  scrambling procedure.
From this point of view the idea underlying Definition~\ref{defasnet} is to preserve the ``good'' dependencies of the points of the $(t,m,d)$-net (namely the ones that result in fairness with respect to elementary intervals of volume $b^{m-t}$) and, in case we have $t>0$, to break up the ``bad'' dependencies (that is dependencies of the $(m-t + 1)$-th, $(m-t+2)$-th, ... digits in each component that do not contribute to the uniformness of the point distribution).

A number of concrete realizations of 
randomized nets that satisfy our definition of an abstract scrambled net have been proposed in the literature, see, e.g., \cite{Owe95,Mat98,LL02, Owe03} and the references mentioned there.
Now we describe the notion of \emph{scrambling of depth~$\ell$} as discussed in 
\cite{Mat98,Mat10}, comment afterwards in Remark~\ref{Rem:History_Scrambling} on other versions of scrambled nets, and demonstrate in Remark~\ref{Rem:Symmetrization} the importance of symmetrization.
In the following,  we essentially blend the definitions and notions used in 
\cite[Sect.~2.4]{Mat10} and  \cite{Owe03}.

For given $b,\ell \in\N$, a mapping $\sigma: [0,1)\to [0,1)$ is called a \emph{$b$-ary scrambling
of depth~$\ell$} if it maps for all $k\in [1..\ell]$ each 
elementary interval $E$ in base $b$ of order $k$ into an elementary
interval $E'$ in the same base having the same order, and distinct $E$ are mapped into distinct $E'$.

The next statement is easily verified, cf. \cite[Proposition~1]{Owe95} or \cite[Observation~2.11]{Mat10}.

\begin{lemma}\label{Lemma:d_tuple_scram}
If $P\subset [0,1)^d$ is a $(t,m,d)$-net in base $b$ and
$\bssigma = (\sigma_1, \ldots, \sigma_d)$ is a $d$-tuple of $b$-ary
scramblings of depth $m-t$, then 
$\bssigma(P)= \{ \bssigma(\bsp) \,|\, \bsp\in P \}$ with 
$\bssigma(\bsp) := (\sigma_1(p_1), \ldots, \sigma_d(p_d))$ is again a $(t,m,d)$-net in base $b$. 
\end{lemma}

The general form of a $b$-ary scrambling $\sigma$ of depth~$\ell$ is 
the following:
Let $x\in [0,1)$ have $b$-ary digit expansion 
\begin{equation}\label{dig_exp}
x= \xi_1 b^{-1} + \xi_2 b^{-2} + \xi_3 b^{-3}+\cdots;
\end{equation}
to ensure that the expansion is unique we assume that infinitely many $\xi_j\in \{ 0,1, \ldots, b-1\}$ are different from $b-1$.
Then
the first $b$-ary digit of $\sigma(x)$ is $\pi(\xi_1)$, where 
$\pi$ is a permutation of $\{0,1,\ldots,b-1\}$. For $j=2,\ldots,\ell$ 
the $j$th digit
of $\sigma(x)$ is given by $\pi_{\xi_1,\ldots,\xi_{j-1}}(\xi_j)$,
where $\pi_{\xi_1,\ldots,\xi_{j-1}}$ is a permutation of 
$\{0,1,\ldots,b-1\}$ that may depend on the first $j-1$ digits of~$x$.
Furthermore, let $\psi: [0,1)\to [0,1)$ be a mapping. Then the general form of a $b$-ary scrambling $\sigma$ of depth~$\ell$ is 
$$
\sigma(x) = \sum_{j=1}^\ell \pi_{\xi_1,\ldots,\xi_{j-1}}(\xi_j) b^{-j}
+ \psi(x) b^{-\ell}.
$$
A \emph{nested uniform $b$-ary random scrambling of depth $\ell$} is defined
in the following way: We pick permutations $\pi_{\xi_1,\ldots,\xi_j}$,
$j=0,1,\ldots, \ell-1$, $\xi_1,\ldots,\xi_{\ell-1} \in \{0,1,\ldots, b-1\}$,
independently and uniformly at random. To map a finite multiset $F=(x^{(n)})_{n=1}^N$
of $[0,1)$,
we define $\sigma(x)$ for each $x=x^{(n)}$ in  $F$ with $b$-ary representation as 
in \eqref{dig_exp} by
\begin{equation}\label{form_nested}
\sigma(x) = \sum_{j=1}^\ell \pi_{\xi_1,\ldots,\xi_{j-1}}(\xi_j) b^{-j}
+ y_x b^{-\ell},
\end{equation}
where $y_x = y_{x^{(n)}} \in [0,1)$ is chosen uniformly at random and independently for  $n=1,\ldots,N$.
We emphazise that $y_{x^{(i)}}$ and $y_{x^{(j)}}$ are independent random variables for $i \neq j$, even if $x^{(i)} = x^{(j)}$.

Nested uniform $b$-ary random scrambling is essentially what was 
proposed by Owen in \cite{Owe95}, see also Remark~\ref{Rem:History_Scrambling}.

A simplified version of nested uniform $b$-ary  random scrambling of depth $\ell$
is \emph{positional uniform $b$-ary random scrambling of depth $\ell$}
proposed by Matou\v{s}ek in \cite{Mat98}: Here one only picks for each ``digit position'' $j\in [\ell]$
a permutation $\pi_j$ of $\{0,1,\ldots,b-1\}$ independently
and uniformly at random and defines for $x = x^{(n)}$ in $F$ as in \eqref{dig_exp}
\begin{equation}\label{form_positional}
\sigma(x) = \sum_{j=1}^\ell \pi_{j}(\xi_j) b^{-j} + y_x b^{-\ell},
\end{equation}  
where, as in \eqref{form_nested}, $y_x = y_{x^{(n)}}$ is chosen uniformly at random and independently for distinct $n\in [1..N]$. 

So in the definitions given above ``random scrambling of depth $\ell$'' refers to the treatment of the different digits: for each $x\in F$  the first $\ell$ digits of $\sigma(x)$ are obtained by applying random permutations to the first $\ell$ digits of $x$, but the remaining digits of $\sigma(x)$ are determined by independent random variable $y_x b^{-\ell}$ in $[0,b^{-\ell})$, cf. \eqref{form_nested} and \eqref{form_positional}.
The attributes ``nested'' and ``positional'' refer to the \emph{scrambling framework of the first $\ell$ digits}: ``nested'' refers to the case where the permutation used to scramble the $j$th digit of  $x\in [0,1)$ as in \eqref{dig_exp} depends (in a nested fashion) on all the previous digits of $x$ and ``positional'' refers to the case where it only depends on the digit position $j$ itself (independently of the previous digits of $x$).
In both cases it is understood that all permutations are chosen independently at random in some subset of the set of all permutations of $\{0,1,\ldots, b-1\}$. 
The attribute ``uniform'' now refers 
the specific choice of the subset of permutations and the chosen probability distribution on it: we choose them uniformly at random from the set of all permutations of $\{0,1,\ldots,b-1\}$. Other choices are possible, as, e.g., the sets of all \emph{linear permutations}, as suggested in \cite{Mat98},  or the set of all \emph{digital shifts}, as mentioned in \cite{LL02}. These sets, endowed with the corresponding uniform probability distribution, yield \emph{nested} and \emph{positional linear} or \emph{nested} and \emph{positional digital shift $b$-ary random scrambling of depth $\ell$}, respectively.
Surely, there exist more forms of random scrambling. One example is \emph{affine matrix $b$-ary random scrambling of depth $\ell$}, i.e., applying \emph{affine Matrix random scrambling} as proposed in \cite{Mat98} to the first $\ell$ digits of $x\in F$ as in \eqref{dig_exp}: for $k=1,\ldots, \ell$ the $k$th digit of $\sigma(x)$ is given by
\begin{equation}\label{mat_scram}
\sigma(x)_k = h_{k,k} \xi_k + \Bigg( \sum_{j=0}^{k-1} h_{k,j} \xi_j + g_k \Bigg) =: \pi^{\rm mat}_{\xi_1, \ldots, \xi_{k-1}}(\xi_k),
\end{equation}
where all $h_{k,k}\in \{1, \ldots, b-1\}$, $h_{k,j} \in \{0,1,\ldots, b-1\}$, $1\le j < k$,
$g_k \in \{0,1,\ldots, b-1\}$ are uniformly distributed and mutually independent.
This is in some sense a hybrid of nested and positional random scrambling: the resulting random permutations $\pi^{\rm mat}_{\xi_1, \ldots, \xi_{k-1}}$ are clearly nested, but 
the randomization is positional, i.e., the random variables $h_{1,k}, \ldots, h_{k,k}$ and $g_k$ in \eqref{mat_scram} that generate the random permutations for the $k$th digits are only depending on the digit position $k$. 

We mention random scramblings that are not ``of depth $\ell$'' in Remark~\ref{Rem:History_Scrambling}. 

\begin{observation}\label{Observation}
If one wants to construct abstract scramblings by applying 
a $d$-tuple ${\bssigma} = (\sigma_1, \ldots, \sigma_d)$ of independent
$b$-ary random scramblings of depth $m-t$ of the form \eqref{form_nested} to a $(t,m,d)$-net $P$ in base $b$, 
one has to fulfill conditions (i), (ii), and (iii) of Definition~\ref{defasnet}.
We make the following observations:
\begin{enumerate}
\item Condition~(i) is satisfied due to Lemma~\ref{Lemma:d_tuple_scram}. 
\item To ensure condition (ii) it is sufficient that for
each $j\in [1..d]$ all random permutations $\pi$ employed by $\sigma_j$ have 
\emph{single-fold uniformity} \cite[Definition~2.5]{Owe03}
$$\PP (\pi(i) = j) = 1/b \hspace{3ex}\text{for all $i,j\in \{0,1,\ldots,b-1\}$,}$$
and all families of these random permutations that scramble different digits are mutually independent.
This is obviously satisfied for nested or positional random scrambling if we choose the permutations uniformly at random from the sets of all permutations, of all linear permutations, or of all digital shifts, respectively (cf. also \cite{Owe03}).
Furthermore, it is satisfied for affine matrix random scrambling (cf. also \cite{Mat98}).
\item Condition (iii) is satisfied, since the position of each $\bssigma(\bsp)$ inside the  basic cube it ends up in is solely determined by the additive vectors $\bsy_{\bsp} b^{t-m}$ that are independent and uniformly distributed in $[0,b^{t-m})^d$,
cf. \eqref{form_nested} and \eqref{form_positional}. 
\end{enumerate}
\end{observation}

We summarize the findings from Observation~\ref{Observation} in the following Proposition.


\begin{proposition}
\label{ElemProp}
Let $P = \{\bsp_1, \dots, \bsp_N\} \subset [0,1)^d$ be a $(t,m,d)$-net in base $b$. Let $\bssigma=(\sigma_1,\ldots,\sigma_d)$ be a $d$-tuple of mutually independent 
random scramblings of depth $m-t$ of the form \eqref{form_nested}, where for each $j\in [1..d]$ all random permutations employed by $\sigma_j$ have single-fold uniformity and all families of random permutations that scramble different digits are mutually independent.
Then $(\bssigma(p_i))_{i=1}^N$ is an abstract scrambled $(t,m,d)$-net.

The statement holds, in particular, if each $\sigma_j$, $j=1,\ldots,d$, is a 
nested or positional (i) uniform or (ii) linear or (iii) digit shift
$b$-ary random scrambling of depth $m-t$, or an affine matrix $b$-ary random scrambling of depth $m-t$.
%
\end{proposition}

\begin{remark}\label{Rem:History_Scrambling}
Crucial for the realization of abstract scramblings as $b$-ary random scramblings of $(t,m,d)$-nets in base $b$ is the concept of scrambling of depth $m-t$.
When introducing $b$-ary random scrambling in \cite{Owe95}, Owen first proposed 
infinite-digit scrambling
\begin{equation}\label{inf_dig_scram}
\sigma(x) = \sum_{j=1}^\infty \pi_{\xi_1,\ldots,\xi_{j-1}}(\xi_j) b^{-j}
\end{equation}
(but for computational issues he also proposed scrambling of depth $\ell$ for $\ell \ge m$, cf. \cite[Section~3.3]{Owe95}).
Assume that all random permutations considered have single-fold uniformity and are mutually independent.
Then it is easy to see that for $(0,m,d)$-nets using $d$-tuples of scramblings of the form \eqref{inf_dig_scram} and \eqref{form_nested} with $\ell \ge m$ lead to the same random scrambling. But this is not true any more if we consider $(t,m,d)$-nets with parameter $t>0$; then \eqref{form_nested} with $\ell = m-t$ yields an abstract scrambling, but \eqref{inf_dig_scram} may violate condition~(iii) of Definition~\ref{defasnet}. 
\end{remark}

\begin{remark}\label{Rem:Symmetrization}
Consider for $m\in \N$ and $N=b^m$ an arbitrary point set $P = (p^{(n)})_{n=1}^N$ in $[0,1)^d$, and
let $\bssigma = (\sigma_1, \ldots, \sigma_d)$ be a $d$-tuple of $b$-ary
scramblings of depth $m-t \ge 1$ or Owen's infinte digit scrambling discussed in Remark~\ref{Rem:History_Scrambling}. Then the random point set
$\bssigma(P)= ( \bssigma(p^{(n)}))_{n=1}^N$ cannot be $\mathcal{C}^d_0$-$\gamma$-negatively dependent for any $\gamma$ independent of $N=b^m$. Indeed, due to the pidgeon hole principle we find an elementary box $E$ of volume $b^{-1}$ that contains
at least $N/b$ points of $P$, without loss of generality the points $p^{(1)}, \ldots, p^{(N/b)}$. Let $c$ be the lower left corner of $E$. Put $E_0 := E-c$, i.e., $E_0$ is the elementary box that has the same shape as $E$ and its lower left corner in $0$. 
Then we have for all $j\in [1.. N]$ that 
$$
\PP  \big(  \bssigma(p^{(j)}) \in E_0 \big)= b^{-1},
$$
and, since $m-t \ge 1$, 
\begin{equation*}
\begin{split}
&\PP \big(  \bssigma(p^{(1)}), \ldots,  \bssigma(p^{(N/b)}) \in E_0 \big) \\
= &\PP  \big(  \bssigma(p^{(1)}) \in E_0 \big) \cdot
\PP \big(  \bssigma(p^{(2)}), \ldots,  \bssigma(p^{(N/b)}) \in E_0 \, \big| \, \bssigma(p^{(1)}) \in E_0 \big) 
= b^{-1} \cdot 1
= b^{-1}.
\end{split}
\end{equation*}
This immediately yields
$$\gamma_{\mathcal{D}^d_0} (\bssigma(P)) \ge \gamma_{\mathcal{C}^d_0} (\bssigma(P)) \ge \frac{1}{b} \cdot b^{N/b},$$
showing that those correlation numbers grow at least exponentially in the number of points $N$.
The unfavorable correlation numbers are due to the dependencies of the first $m-t$ digits of all components of the original points that are preserved by the specific randomizations we consider. 
The remedy to this problem is simple -- we just have to symmetrize the random points, i.e., we have to ensure that they are exchangeable random variables. This can be done by simply giving each point of $P$ a random label at the beginning: Choose a permutation $\pi$ uniformly at random and put $X_j:= \bssigma(p^{(\pi(j))})$ for all $j=1,\ldots,N$.
Now the random point set $\X = (X_j)_{j=1}^N$ contains exactly the same random points as $\bssigma(P)$, but the random points $\bssigma(p^{(\pi(1))}), \ldots, \bssigma(p^{(\pi(N))})$ are now exchangeable. 
In \cite{Lem18, WLD20} the authors do random scrambling without choosing random labels; instead the symmetrization procedure is  implicitly contained in the definition of the negative dependence properties studied by those authors, cf. \cite[Sect.~2.1]{Lem18} and  \cite[Sect.~2.2]{WLD20}.\\
We emphasize that our lower bound for the correlation number in Proposition~3.9 holds for abstract scrambled nets, regardless if they are symmetrized or not.
\end{remark}


Since every Latin hypercube sample with $N$ points is, in particular, an abstract scrambled $(0,1,d)$-net in base $b:=N$,  we see from Proposition~\ref{prop33} that there does not exist a constant $\gamma \in [1,\infty)$ such that all $(0,m,d)$-nets in base $b$ are $\mathcal{D}^d_0$-$\gamma$-negatively dependent. In particular, not all 
scrambled $(0,m,d)$-nets are $\mathcal{D}^d_0$-$1$-negatively dependent.

The next example shows that there exist scrambled $(0,m,d)$-nets in base $b$ different from Latin hypercube samples that are not $\mathcal{D}^d_0$-$1$-negatively dependent.

\begin{proposition}\label{Ex:Correlation_Constant_Nets}
Consider an abstract scrambled $(0,2,3)$-net $\Y = (Y_i)^N_{i=1}$ 
in base~$2$. 
Then for all $\varepsilon > 0$ 
there is a $D \in \mathcal{D}^d_0$ such that
\begin{equation*}\label{quotient_2}
\frac{\PP(Y_1,Y_2, Y_3\in D)}{\prod_{i=1}^3 \PP(Y_i\in D)} \ge \left( \frac{4}{3} \right)^3 - \varepsilon = 2.\overline{370} - \varepsilon.
\end{equation*}
In particular, $\gamma_{\mathcal{D}^3_0}(\Y) \ge 2.\overline{370}$.
\end{proposition}

Recall that essentially all commonly used ways of generating randomly scrambled $(0,m,d)$-nets, including Owen's originally proposed infinite digit-scrambling and scrambling of depth $\ell$ for $\ell \ge m$, result in an abstract scrambled $(0,m,d)$-net, cf. Proposition~\ref{ElemProp} and Remark~\ref{Rem:Symmetrization}.
 

\begin{proof}
  We note first that it follows right from the definitions that each $(0,m,d)$ net $P$ in base $b$ inherits the property from Latin hypercube samples that for every $j \in [1..d]$ and $k \in [1..N]$ there is exactly one point of $P$ in each box $E_{jk} = \prod_{i=1}^d E_i$ with $E_j = [\frac{k-1}{N},\frac{k}{N})$ and $E_i = [0,1)$ for $i \neq j$. In particular, $P$ cannot contain two distinct points $Y$ and $Y'$ such that there are $j \in [1..d]$ and $k \in [1..N]$ with $Y_j, Y'_j \in [\frac{k-1}{N},\frac{k}{N})$. We shall use this elementary \emph{combinatorial property} frequently in the following.

Let now $(Y_1, \dots, Y_4)$ be an abstract scrambled $(0,2,3)$ net in base~$2$. 
That is, we have $N=4$ points in $[0,1)^d$, where $d=3$, such that each elementary interval in base~$2$  of order~$2$ contains exactly one point.
Let  $\varepsilon \in (0,\frac{1}{4})$.
\begin{equation*}
A:= [0,\tfrac{1}{2})^3,  B:= [0, \tfrac{1}{2} + \varepsilon)^3, 
\hspace{1ex}\text{and}\hspace{1ex}
D:= B\setminus A.
\end{equation*}
Denote by $\EE$ the event that each of $Y_1$, $Y_2$, and $Y_3$ has exactly one coordinate entry in $[\frac{1}{2}, \frac{3}{4})$ and the other two coordinate entries in $[0, \frac 12)$. Furthermore, put $U:= [\frac{3}{4}, 1)^3$. 
We first show that the combinatorial property implies the following inclusions of events:
\begin{equation}
\{Y_1, Y_2, Y_3 \in D\} \subseteq \EE \subseteq \{Y_4 \in U\}.  \label{eq:inclusions}
\end{equation}
Indeed, assume that $Y_1,Y_2, Y_3 \in D$. Then each of these three points has at least one coordinate  entry in $[\frac{1}{2}, \frac{3}{4})$ and all coordinate entries not in $[\frac 12, \frac 34)$ are in $[0,\frac 12)$. 
By the combinatorial property established in the beginning of this proof, each point can have at most one coordinate in $[\frac 12, \frac 34)$. Hence each point $Y_1$, $Y_2$, and $Y_3$ has exactly one coordinate entry in $[\frac{1}{2}, \frac{3}{4})$, establishing the first inclusion.

To prove the second one, assume that $\EE$ holds. By the combinatorial property and elementary counting, we see that for each $j \in [1..3]$ and $k \in [1..3]$, there is a point $Y_i$, $i \in [1..3]$, such that $(Y_i)_j \in [\frac {k-1}{4}, \frac k4)$. Hence $Y_4$ cannot have any coordinate in $[0,\frac 34)$, showing that $Y_4 \in U$. 
%


The crucial observation, following from the stronger net property and leading to an easy computation of $\PP (Y_1,Y_2,Y_3 \in D \,|\, Y_4 \in U)$, is that the second inclusion in~\eqref{eq:inclusions} is in fact an equality, that is, we have 
\begin{equation}
\EE = \{Y_4 \in U\}. \label{eq:equality}
\end{equation} 
To show~\eqref{eq:equality}, assume that $Y_4 \in U$. Since the canonical box $B = [\frac 12,1) \times [\frac 12,1) \times [0,1)$ contains $Y_4$, none of the other three points can be contained in $B$. In particular, none of these three points can have both the first two coordinate entries in $[\frac 12, \frac 34)$. The same argument applied to the canonical boxes $[\frac 12,1) \times [0,1) \times [\frac 12,1)$ and $[0,1) \times [\frac 12,1) \times [\frac 12,1)$ shows that each of $Y_1, Y_2, Y_3$ can have at most one coordinate entry in $[\frac 12, \frac 34)$, and by the combinatorial property this becomes exactly one coordinate entry in $[\frac 12, \frac 34)$. Again by the combinatorial property and $Y_4 \in U$, the other two coordinate entries of each of these points must lie in $[0,\frac 12)$. This establishes $\EE$ and thus~\eqref{eq:equality}.

Let $\CC$ be the set of all $C = (C_1, \dots, C_4)$ such that (i)~$C_1, \dots, C_4$ are basic cubes, (ii)~$C_4 = U$, and (iii)~the event $\EE_C = \{\forall i \in [1..4] : Y_i \in C_i\}$ has positive probability. By~\eqref{eq:equality}, the event $\{Y_4 \in U\}$ implies the event $\EE$. Hence also any event $\EE_C$, $C \in \CC$ implies $\EE$, and thus 
\begin{equation}
  \frac{\lambda^d(D \cap C_i)}{\lambda^d(C_i)} = 4\varepsilon. \label{eq:4eps}
\end{equation} 

By~\eqref{eq:inclusions}, the definition of abstract scrambled nets, and~\eqref{eq:4eps}, we have
\begin{align*}
  \PP(Y_1, Y_2, Y_3 \in D) & = \PP(\{Y_1, Y_2, Y_3 \in D\} \cap \{Y_4 \in U\}) \\
  & = \sum_{(C_1, \dots, C_4) \in \CC} \PP(\{Y_1, Y_2, Y_3 \in D\} \cap \EE_C) \\ 
  & = \sum_{(C_1, \dots, C_4) \in \CC} \PP(Y_1, Y_2, Y_3 \in D \mid \EE_C) \PP(\EE_C) \\ 
  & = \sum_{(C_1, \dots, C_4) \in \CC} \prod_{i=1}^3 \frac{\lambda^d(D \cap C_i)}{\lambda^d(C_i)} \cdot \PP(\EE_C) \\ 
  & = \sum_{(C_1, \dots, C_4) \in \CC} (4\varepsilon)^{3} \PP(\EE_C) \\ 
  & = (4\varepsilon)^{3} \PP(Y_4 \in U) = (4\varepsilon)^{3} \tfrac 1 {64} = \varepsilon^3,
\end{align*}
  where in the last line we obtained $\PP(Y_4\in U) = \frac{1}{64}$ from the fact that every point in an abstract scrambled net is uniformly distributed in $[0,1)^d$.

We easily compute 
$$\lambda^3(D) = \left( \frac{1}{2} + \varepsilon \right)^3 -  \left( \frac{1}{2} \right)^3 = \varepsilon \left( \frac{3}{4} + \frac{3}{2} \varepsilon + \varepsilon^2 \right) = \frac{3}{4} \varepsilon + \Theta(\varepsilon^2).$$
We thus have $\PP(Y_i \in D) = \lambda^3(D) = \varepsilon \left(\frac 34 + \Theta(\varepsilon)\right)$ for all $i \in [1..3]$. 
Consequently, for $\varepsilon \searrow 0$ we have 
\begin{equation*}
\frac{\PP(Y_1,Y_2, Y_3\in D)}{\prod_{i=1}^3 \PP(Y_i\in D)} = \left( \frac{4}{3} \right)^3 \cdot \left( \frac{1}{1+\Theta(\varepsilon)} \right)^3 \nearrow  \left( \frac{4}{3} \right)^3 = 2.\overline{370}.
\end{equation*}
\end{proof}

\begin{remark}
As a comparison, we want to study the setting that leads to the unfavorable lower bound of the correlation constant for scrambled nets in Proposition~\ref{Ex:Correlation_Constant_Nets} now for Latin hypercube samples $(X_i)^{4}_{i=1}$ in $[0,1)^3$. We will observe that in this setting indeed the stricter dependencies of the scrambled net lead to a larger lower bound on the $\mathcal{D}^3_0$-correlation number. We use the same notation as in Proposition~\ref{Ex:Correlation_Constant_Nets} and the corresponding proof.

 As in~\eqref{eq:inclusions}, we have
 \begin{equation*}
\{X_1, X_2, X_3 \in D\} \subseteq \EE \subseteq \{X_4 \in U\}.
\end{equation*}
This gives us, in particular,
\begin{align*}
\PP (X_1,X_2,X_3 \in D) & = \PP (X_1,X_2,X_3 \in D \,|\, X_4 \in U) \cdot \PP(X_4\in U). 
\end{align*}
We also have again $\PP(X_4\in U) = \frac{1}{64}$. What is different is that $X_4 \in U$ does not imply~$\EE$. In fact, if we have $X_4 \in U$, then, apart from scaling, $X_1, X_2, X_3$ form a Latin hypercube sample in $[0,\frac 34)^3$, and as computed in the proof of Proposition~\ref{prop33}, the event $\EE$ occurs with probability $\frac{3!}{3^3} = \frac 29$. Consequently,
\begin{align*}
\PP (X_1,X_2,X_3 \in D \,|\, X_4 \in U) &= \PP (X_1,X_2,X_3 \in D \,|\, \EE) \cdot \PP(\EE \,|\, X_4 \in U) \\
&= \tfrac 29 \PP (X_1,X_2,X_3 \in D \,|\, \EE) = \tfrac 29 (4\varepsilon)^3.
\end{align*}
We thus have $\PP (X_1,X_2,X_3 \in D) = \frac 29 \varepsilon^3$.
If we let $\varepsilon$ tend to zero, we obtain
\begin{equation*}
\frac{\PP(X_1,X_2, X_3\in D)}{\prod_{i=1}^3 \PP(X_i\in D)} =  \frac{2}{9} \cdot \left( \frac{4}{3} \right)^3 \cdot \left( \frac{1}{1+\Theta(\varepsilon)} \right)^3 \nearrow   \frac{2}{9} \cdot \left( \frac{4}{3} \right)^3 
=  0.5267\ldots.
\end{equation*}
\end{remark}

\subsection*{Acknowledgment}
This work was initiated at the Dagstuhl Seminar 19341, ``Algorithms and Complexity for Continuous Problems''. Part of this work was done while Michael Gnewuch was a ``Chercheur Invit\'e'' of \'Ecole Polytechnique. He thanks his colleagues at the LIX for their hospitality.

\footnotesize


\end{document}